\documentclass[reqno, 11pt]{amsart}
\usepackage{amssymb,verbatim}
\usepackage[nesting]{hyperref}

\usepackage[pdftex]{graphicx}
\usepackage{listings}
\usepackage{multirow}
\usepackage{placeins}
\usepackage{color}
\usepackage{subfigure}
\usepackage{lscape}
\usepackage{dsfont}
\usepackage{bbm}

\usepackage[T1]{fontenc}
\usepackage[latin9]{inputenc}
\usepackage{xcolor}
\usepackage{pdfcolmk}

\usepackage{units}
\usepackage{mathtools}
\usepackage{amsthm}
\usepackage{amsmath}
\usepackage{amssymb}
\PassOptionsToPackage{normalem}{ulem}
\usepackage{ulem}

\makeatletter

\newcommand{\BlackBoxes}{\global\overfullrule5pt}

\BlackBoxes

\newcommand{\R}{\mathbb{R}}
\newcommand{\C}{\mathbb{C}}

\newcommand{\N}{\mathbb{N}}

\newcommand{\cH}{\mathcal{H}}
\newcommand{\cS}{\mathcal{N}}
\newcommand{\cN}{\mathcal{N}}

\newcommand{\cE}{\mathcal{E}}
\newcommand{\cL}{\mathcal{L}}
\newcommand{\cU}{\mathcal{U}}
\newcommand{\cV}{\mathcal{V}}
\newcommand{\cP}{\mathcal{P}}

\numberwithin{equation}{section}
\numberwithin{figure}{section}
\theoremstyle{plain}
\newtheorem{thm}{\protect\theoremname}
\theoremstyle{plain}
\newtheorem{prop}[thm]{\protect\propositionname}
\theoremstyle{definition}
\newtheorem{defn}[thm]{\protect\definitionname}
\theoremstyle{remark}
\newtheorem{rmk}[thm]{\protect\remarkname}
\theoremstyle{plain}
\newtheorem{lem}[thm]{\protect\lemmaname}
\theoremstyle{plain}
\newtheorem{cor}[thm]{\protect\corollaryname}
\newtheorem{ex}[thm]{\protect\examplename}

\newcommand{\pf}{\noindent \textit{Proof}:\ }
\usepackage{datetime,color,currfile,pdfsync,cancel,hyperref}

\providecommand{\corollaryname}{Corollary}
\providecommand{\definitionname}{Definition}
\providecommand{\lemmaname}{Lemma}
\providecommand{\propositionname}{Proposition}
\providecommand{\remarkname}{Remark}
\providecommand{\theoremname}{Theorem}
\providecommand{\examplename}{Example}

\def\0{\kern0pt\-\nobreak\hskip0pt\relax}

\makeatletter
\AtBeginDocument{%
 \def\@serieslogo{%
 \vbox to\headheight{%
 \parindent\z@ \fontsize{6}{7\p@}\selectfont
 \vss}}}

\def\makeoverbar#1#2#3#4#5#6#7{%
 \setbox0=\hbox{$\m@th#2\mkern#5mu{{}#3{}}\mkern#6mu$}%
 \setbox1=\null \dimen@=#4\fontdimen8#13 \dimen@=3.5\dimen@
 \advance\dimen@ by \ht0 \dimen@=-#7\dimen@ \advance\dimen@ by \wd0
 \ht1=\ht0 \dp1=\dp0 \wd1=\dimen@
 \dimen@=\fontdimen8#13 \fontdimen8#13=#4\fontdimen8#13
 \rlap{\hbox to \wd0{$\m@th\hss#2{\overline{\box1}}\mkern#5mu$}}
 \fontdimen8#13=\dimen@}

\def\mylabel#1#2{{\def\@currentlabel{#2}\label{#1}}}

\begin{document}



\title[Characterization of symbols of operators]{An analytic characterization of symbols of operators on non-Gaussian Mittag-Leffler functionals}

\author[W. \smash{Bock}]{Wolfgang Bock${}^1$}
\address[W. Bock]{Linnaeus University, Department of Mathematics, Universitetsplatsen 1, 352 52 V\"axj\"o, Sweden}

\email{\href{mailto:wolfgang.bock@lnu.se}
{wolfgang.bock@lnu.se}}

\author[A. \smash{Gumanoy}]{Ang Elyn Gumanoy${}^2$}
\address[A. Gumanoy]{Department of Mathematics, MSUIIT, Iligan, The Philippines}

\email{\href{mailto:angelyn.gumanoy@g.msuiit.edu.ph} {angelyn.gumanoy@g.msuiit.edu.ph}}

\author[S. \smash{Menchavez}]{Sheila Menchavez${}^3$}
\address[S. Menchavez]{Department of Mathematics, MSUIIT, Iligan, The Philippines}

\email{\href{mailto:sheila.menchavez@g.msuiit.edu.ph} {sheila.menchavez@g.msuiit.edu.ph}}

\author[E. \smash{Nabizadeh}]{Elmira Nabizadeh Morsalfard${}^4$}
\address[E. Nabizadeh]{Linnaeus University, Department of Mathematics, Universitetsplatsen 1, 352 52 V\"axj\"o, Sweden}

\email{\href{mailto:elmira.nabizadehmorsalfard@lnu.se}
{elmira.nabizadehmorsalfard@lnu.se}}

\thanks{${}^{1,4}$ Linnaeus University, Department of Mathematics, Universitetsplatsen 1, 352 52 V\"axj\"o, Sweden}
\thanks{${}^{2,3}$Department of Mathematics, MSU-IIT, Iligan, The Philippines}

\begin{abstract}
In this paper, we provide proofs for the analytic characterization theorems of the operator symbols utilizing the characterization theorem for the Mittag-Leffler distribution space. We work out examples which can be interpreted as integral kernel operators and treat the important case of the translation operator.
\end{abstract}
\maketitle

\vspace{0.5cm}
\begin{minipage}{14cm}
{\small
\begin{description}
\item[\rm \textsc{ Key words} Mittag-Leffler measure, operator symbol, $S_{\mu_{\beta}}$-transform, \\
characterization theorem
]
{\small }
\end{description}
}
\end{minipage}

\section{Introduction}
In Gaussian analysis and especially white noise analysis, various characterization theorems were proven to provide a deep understanding of the structure of test and distribution spaces \cite{PS91, KLPSW96, GKS97}. The main ingredient is an infinite dimensional analogue of the Laplace transform called the $S$-transform. The $S$-transform establishes a one-to-one correspondence between the distribution space and the spaces of holomorphic functions on a locally convex space satisfying a certain growth condition. With this theoretical background, many applications in applied mathematics, stochastic analysis, stochastic (partial) differential equations, mathematical physics and many others can be treated, see e.g.~the monographs \cite{HKPS93, Ob94, Kuo96,LetUsUse}. 
The characterization theorems have been generalized to linear bounded operators on the corresponding spaces, see e.g. \cite{Ob94} and \cite{Ch97}, where the symbol of operators was introduced as a generalization of the S-transform. Also we refer to the works of Hida et al. \cite{HOS92}, Un Cig Ji et al \cite{CJ97, CJO02} and Ji, Obata and Ouerdiane \cite{JOO02, OO11}  which treat speacial classes of (pseudo-)differential operators in the Fock space.\\
Mittag-Leffler analysis, initiated by Grothaus et. al. \cite{GJRdS15}, has been established as an infinite dimensional analysis with respect to non-Gaussian measures of Mittag-Leffler type referred to as Mittag-Leffler measures. Using the Bochner and Minlos theorem, such measures were defined using their characteristic function represented by the Mittag-Leffler function--a generalization of the exponential function.  The grey noise measure \cite{Schneider90, Sch92, MM09} is included as a special case in the class of Mittag-Leffler measures, which gives rise to various applications of fractional diffusion equations, which carry numerous applications in science, like relaxation type differential equations or viscoelasticity.  Special classes of operators and their relation to multiplication with special Mittag-Leffler-distributions \cite{BG21}. \\
In this short note, we want to study integral kernel operators and the symbols of general linear operators in the Mittag-Leffler setting. We give a complete characterization based on the characterization theorem. We work out examples which can be interpreted as integral kernel operators and treat the important case of the translation operator.

\section{Preliminaries}
\subsection{Holomorphy on locally convex spaces}
{\rm{In this section, we recall concepts pertaining to holomorphic functions in locally convex topological vector spaces $\cE$ (over the complex field $\C$). These definitions will be beneficial for subsequent discussions \cite{Bar85}, \cite{Col82} and \cite{Din81}.

Let $\cL(\cE^n)$ be the space of $n$-linear mappings from $\cE^n$ into $\C$ and $\cL_s(\cE^n)$ be the subspace of symmetric $n$-linear forms. Moreover, let $P^n(\cE)$ be the set of all $n$-homogeneous polynomials on $\cE$. There is a linear bijection
$$\cL_s(\cE^n) \ni A \longleftrightarrow \hat{A} \in P^n(\cE).$$ }}
\defn{Let $\cU \subset \cE$ and $f:\cU\rightarrow \C$ a function. Then $f$ is said to be $G$-\textit{holomorphic} (or \textit{G$\hat{a}$teaux-holomorphic}) if and only if for all $\theta_0 \in \cU$ and for all $\theta \in \cE$, the mapping
$$\C \ni \lambda \mapsto f(\theta_0 + \lambda \theta) \in \C,$$
is holomorphic in some neighborhood of zero in $\C$.}

\rmk{If $f$ is $G$-holomorphic, then there exists for every $\eta \in \cU$ a sequence of homogeneous polynomials $\frac{1}{n!}d^n f(\eta)$ such that 
$$f(\theta +\eta)=\sum_{n=0}^{\infty}\frac{1}{n!}\widehat{d^n f(\eta)}(\theta),$$
for all $\theta$ from some open neighborhood $\cV$ of zero where $\widehat{d^n f(\eta)}(\theta)$ is the $n$-th partial derivative of $f$ at $\eta$ in direction $\theta$. }

\defn{Let $f:\cU \rightarrow \C$ be a $G$-holomorphic function. Then
\begin{enumerate}
\item[i)] $f$ is said to be \textit{holomorphic}, if and only if for all $\eta \in \cU$ there exists an open neighborhood $\cV$ of zero such that
$$\cV \ni \theta \mapsto \sum_{n=0}^{\infty}\frac{1}{n!}\widehat{d^n f(\eta)}(\theta),$$
converges uniformly on $\cV$ to a continuous function.
\item[ii)] $f$ is \textit{holomorphic} at $\theta_0$ if and only if there is an open set $\cU$ containing $\theta_0$ such that $f$ is holomorphic on $\cU$.
\item[iii)] $f$ is called \textit{entire} if and only if $f$ is holomorphic on $\cE$.
\end{enumerate}}

\prop{{\rm\cite{Din81}}\label{hol-lbbd}
Let $\cU \subset \cN_{\C}$ be open and $f:\cU \rightarrow \C$. Then, $f$ is holomorphic if and only if it is $G$-holomorphic and locally bounded, i.e., each point $\xi\in \cU$ has a neighborhood whose image under f is bounded.
}

\rm{We denote by $\text{Hol}_{0}(\cN_{\C})$ the space of all holomorphic functions at zero. Let $f$ and $\tilde{f}$ be holomorphic on a neighborhood $\cV,\cU \subset \cN_{\C}$ of zero, respectively. We identify $f$ and $\tilde{f}$ if there exists a neighborhood $\mathcal{W}\subset \cV$ and $\mathcal{W}\subset \cU$ such that $f(\xi)= \tilde{f}(\xi)$ for all $\xi \in \mathcal{W}$. $\text{Hol}_{0}(\cN_{\C})$ is the union of the spaces
$$\left\{ f\in \text{Hol}_{0}(\cN_{\C}) \big\vert \text{ n}_{p,l,\infty}(f):=\displaystyle{\sup_{|\xi|_p \leq 2^{-l}}}|f(\xi)|< \infty \right\}, \quad p,l\in \N,$$
and carries the inductive limit topology.

The following corollary is a direct consequence of Proposition \ref{hol-lbbd}.}

\cor{Let $f:\cN_{\C}\rightarrow \C$ be given. Then $f\in {\rm{Hol}}_{0}(\cN_{\C})$ if and only if there exists $p \in \mathbb{N}$, $\varepsilon >0$, and $C > 0$ such that
\begin{enumerate}
\item[i)] ($G$-holomorphic) for all $\xi_0 \in \cN_{\C}$ with $|\xi_{0}|_p \leq \varepsilon$ and for all $\xi\in \cN_{\C}$, the function of one complex variable $\C \ni \lambda \mapsto f(\xi_{0}+\lambda\xi)\in \C$ is holomorphic at zero, and
\item[ii)] (locally bounded) for all $\xi\in \cN_{\C}$ with $|\xi|_p\leq \varepsilon$, we have $|f(\xi)|\leq C$.
\end{enumerate}}

\section{Mittag-Leffler Analysis}

\subsection{Mittag-Leffler measure}

Given a nuclear triple $\cS\subset \cH \subset \cS'$, the space $\cS'$ is equipped with the $\sigma$-algebra $\mathcal{C}_{\sigma}(\cS')$ generated by the cylinder sets
$$\left\{ \omega\in \cS': \left( \left\langle \omega, \varphi_1\right\rangle, \left\langle \omega, \varphi_2\right\rangle,\dots,\left\langle \omega, \varphi_n\right\rangle \right)\in A\right\},$$
where $n\in \mathbb{N}$, $\varphi_1,\dots,\varphi_n\in \cS$ and $A\in \mathcal{B}(\mathbb{R}^n)$, the Borel $\sigma$-algebra over $\mathbb{R}^n$. Since $\cS$ is a projective limit of a countable number of Hilbert spaces, $\mathcal{C}_{\sigma}(\cS')$ coincides with the Borel $\sigma$-algebra  generated by the weak and strong topologies on $\cS'$, refer \cite{BK95}.

The definition of the Mittag-Leffler measure on $(\cS',\mathcal{C}_{\sigma}(\cS'))$ relies on the Mittag-Leffler function which was introduced as a generalization of the exponential function by G$\rm\ddot{o}$sta Mittag-Leffler in \cite{ML05}, refer also \cite{Wim05a, Wim05b}. This function has been further studied and its generalization was first presented in \cite{Wim05a}.

\begin{defn}\label{MLfunc}\cite{GJRdS15}
 For $0<\beta<\infty$, the \emph{Mittag-Leffer function} $E_{\beta}$ is an entire function defined by its power series
$$E_{\beta}(z):=\sum_{n=0}^{\infty}\frac{z^n}{\Gamma(\beta n+
1)}, \quad z\in \mathbb{C},$$
where $\Gamma$ is the Gamma function. In addition, for $0<\gamma<\infty$, the \emph{generalized Mittag-Leffler function} $E_{\beta,\gamma}$ is an entire function defined by the power series
$$E_{\beta,\gamma}(z):=\sum_{n=0}^{\infty}\frac{z^n}{\Gamma(\beta n+\gamma)}, \quad z\in \mathbb{C}.$$
\end{defn}

\rmk{
Here, $\Gamma$ denotes the \emph{gamma function} defined by
$$\Gamma(z)=\int_0^{\infty} x^{z-1}e^{-x}\text{ d$x$}, \quad \text{ for all $z\in \mathbb{C}$,}$$
which is an extension of the factorial to complex numbers, such that $\Gamma(n+1)=n!$ for $n\in \mathbb{N}$. This shows that $E_1(z)=e^z$ for all $z\in \mathbb{C}$. Moreover, the $\Gamma$-function fulfils $\Gamma(z+1)=z\Gamma(z)$ for all $z\in \C$.}

\lem{ \cite{GJRdS15} For the derivative of the Mittag-Leffler function $E_{\beta}$, $0<\beta<\infty$, it holds
\begin{equation}\label{mittag-derivative}
    \frac{\rm{d}}{\rm{d}z}E_{\beta}(z)=\frac{E_{\beta, \beta}(z)}{\beta}, \quad z\in \C.
\end{equation}}

\rm{We also consider a special case of the family of Wright functions the so-called $M$-Wright function $M_{\beta}$ for $0<\beta\le1$ (in one variable) where its series expansion is given by 
\[
M_{\beta}(z):=\sum_{n=0}^{\infty}\frac{(-z)^{n}}{n!\Gamma(-\beta n+1-\beta)}.
\]
The density $M_{\beta}$ if proper symmetrized generalizes the Gaussian density, for more details see for example \cite{MMP10}. In particular for $\beta=\frac{1}{2}$, the corresponding $M$-Wright function reduces to the Gaussian density
\begin{equation}
M_{\frac{1}{2}}(z)=\frac{1}{\sqrt{\pi}}\exp\left(-\frac{z^{2}}{4}\right).\label{eq:MWright_Gaussian}
\end{equation}

The integral representation for the Mittag-Leffler function is given by 
\begin{equation}
E_{\beta}(-s)=\int_{0}^{\infty} e^{-\tau s}{\rm{d}}\nu_{\beta}(\tau)
\end{equation}
where $\nu_{\beta}$ is a probability measure on $[0,\infty)$ with density $M_{\beta}$ with respect to the Lebesgue measure, that is, $\rm{d}\nu_{\beta}(\tau)= M_{\beta}(\tau)\rm{d}(\tau)$.

The Mittag-Leffler function $E_{\beta}$ and the $M$-Wright are related through the Laplace transform
\begin{equation}
\int_{0}^{\infty}e^{-s\tau}M_{\beta}(\tau)\,d\tau=E_{\beta}(-s).\label{eq:LaplaceT_MWf}
\end{equation}}

\rmk{\cite{Jah15} For $0<\beta\leq 1$, the map $[0,\infty)\ni x\rightarrow E_{\beta}(-x)\in \mathbb{R}$ is completely monotonic, that is,
$$(-1)^nE_{\beta}^{(n)}(-x)\geq 0 \text{ for all $x\geq 0$.}$$
Using this fact and \cite{Pol48}, one can show in a similar manner as that of \cite{Sch92} that the map
$$\cS\ni \varphi \mapsto E_{\beta}\left(-\frac{1}{2}\left\langle \varphi,\varphi\right\rangle\right)\in \mathbb{R}$$
is a characteristic function on $\cS$.}

Using the Bochner-Minlos theorem, see \cite{BK95}, we obtain the following definition.

\begin{defn}\label{MLmeas}
\cite{Jah15} For $0<\beta\leq 1$, the \emph{Mittag-Leffler measure} $\mu_{\beta}$ is defined as the unique probability measure on the space $(\cS',\mathcal{C}_{\sigma}(\cS'))$ whose characteristic function is 
\begin{equation}
\int_{\cS'}e^{i\langle \omega,\varphi\rangle}\,d\mu_{\beta}(\omega)=E_{\beta}\left(-\frac{1}{2}\left\langle \varphi,\varphi\right\rangle\right),\quad\varphi\in \cS.\label{eq:ch-fc-gnm}
\end{equation}
The corresponding $L^p$ spaces of complex-valued functions are denoted by $L^p(\mu_{\beta}):=L^p(\cS', \mu_{\beta};\mathbb{C})$ for $p\geq 1$ with corresponding norm $||\cdot||_{L^p(\mu_{\beta})}$. For $p=2$, the corresponding scalar product is denoted by $((\cdot,\cdot))_{L^2(\mu_{\beta})}$. That is, 
\[
(\!(F,G)\!)_{L^{2}(\mu_{\beta})}:=\int_{\cS'}F(\omega)\overline{G}(\omega)\,d\mu_{\beta}(\omega),\quad F,G\in L^{2}(\mu_{\beta}).
\]
\end{defn}

\rmk{\cite{Jah15} The class of Mittag-Leffler measures on $\cS'$ includes the following: 
\begin{enumerate}
\item For $\beta=1$, the Mittag-Leffler measure $\mu_1$ on $\cS'$ is the usual Gaussian measure on $\mathcal{N}'$ with covariance given by the scalar product in $\mathcal{H}$.
\item The measure $\mu_{\beta}$ on $\cS'(\mathbb{R})$ is called grey noise (reference) measure, cf.\ \cite{GJRdS15} and \cite{GJ16}.
\end{enumerate}
}



\rm{Many properties of the Mittag-Leffler measure $\mu_{\beta}$ follows from \eqref{eq:ch-fc-gnm}. The moments of $\mu_{\beta}$ are given by }
\prop{\label{mlaprop}{\rm\cite{GJRdS15}}
For any $\varphi\in \cS$ and $n\in \mathbb{N}_0$,
\begin{align*}
\int_{\cS'} \left\langle \omega,\varphi \right\rangle^{2n+1} \,d\mu_{\beta}(\omega) &=0, \\
\int_{\cS'}\langle \omega,\varphi\rangle^{2n}\,d\mu_{\beta}(\omega) &=\frac{(2n)!}{2^{n}\Gamma(\beta n+1)}\langle \varphi,\varphi\rangle^n .
\end{align*}
In particular for all $\varphi,\psi\in \cS$ we obtain
\begin{align*}
\mathbb{E}_{\mu_{\beta}}[\langle \cdot, \varphi \rangle^2]=||\left\langle \cdot, \varphi\right\rangle||^2_{L^2(\mu_{\beta})}&=\frac{1}{\Gamma(\beta +1)}|\varphi|_0^2,\\   
\int_{\cS'}\langle \omega,\varphi\rangle\langle \omega,\psi\rangle\,d\mu_{\beta}(\omega) &=\frac{1}{\Gamma(\beta +1)}\langle \varphi,\psi\rangle.
\end{align*}}


\subsection{The biorthogonal system}
In non-Gaussian analysis, a different approach is taken, utilizing what is known as the Appell systems. These are biorthogonal systems consisting of Appell polynomials and the corresponding distributions. We repeat a brief construction of the biorthogonal system with respect to the Mittag-Leffler measure $\mu_{\beta}$. The detailed construction can be found in \cite{GJRdS15} or \cite{KSWY98}.

We define the $\mu_{\beta}$-exponential by
$$e^{\mu_{\beta}}: \cS_{\C} \times \cS'_{\C}\longrightarrow \C, \quad \left(\omega, \theta \right) \mapsto e_{\mu_{\beta}}\left(\omega; \theta \right):= \frac{e^{\langle\omega, \theta \rangle}}{l_{\mu_{\beta}}(\theta)}.$$
It is well-define if and only if $l_{\mu_{\beta}}(\theta)\neq 0$. Since $l_{\mu_{\beta}}(0)=1$, $l_{\mu_{\beta}}$ is holomorphic, this implies that there exists a neighborhood $\cU_0\subset \cS_{\C}$ of zero such that $l_{\mu_{\beta}}(\theta)\neq 0$, for all $\theta\in \cU_0$. 

\rmk{With $l_{\mu_{\beta}}(\theta)=E_\beta\left(\frac{1}{2}\langle \theta, \theta\rangle\right)>0$, the $\mu_{\beta}$-exponential is given by
\begin{equation}\label{ser-exp}
e_{\mu_{\beta}}\left(\omega; \theta \right)= \frac{e^{\langle\omega, \theta \rangle}}{E_\beta\left(\frac{1}{2}\langle \theta, \theta\rangle\right)}.
\end{equation}
We see that $E_{\beta}$ is holomorphic and $E_{\beta}(0)=1$. Thus there exists $\varepsilon_{\beta}>0$ such that $E_{\beta}(z)>0$ for all $z\in \mathbb{C}$, $|z|\leq \varepsilon_{\beta}$. Hence, the function $\mathcal{E}:=\frac{1}{E_{\beta}}$ is holomorphic on $B_{\varepsilon_{\beta}}(0)$. Assume that $\mathcal{E}$ has the series expansion
$$\mathcal{E}(z)=\sum_{n=0}^{\infty}b_nz^n, \quad |z|<\varepsilon_{\beta},$$
for some coefficients $b_n\in \mathbb{C}$ and $\theta \in \mathcal{U}_{\beta}$ where
$$\mathcal{U}_{\beta}:=\left\{ \theta \in \cS_{\C}\middle| \frac{1}{2}|\langle \theta, \theta\rangle| <\varepsilon_{\beta}\right\} \subset \cS_{\C}.$$
Then, 
\begin{equation}\label{ser-hol}
\mathcal{E}\left( \frac{1}{2}\langle \theta, \theta\rangle\right) = \sum_{n=0}^{\infty} \frac{b_n}{2^n}\langle \tau, \theta^{\otimes 2}\rangle^n =  \sum_{n=0}^{\infty} \frac{b_n}{2^n}\langle \tau^{\otimes n}, \theta^{\otimes 2n}\rangle. 
\end{equation}
Here $\tau\in (\cS_{\C}^{\otimes 2})'$ denotes the trace operator  which is uniquely defined by 
$$\langle \tau, \theta \otimes \eta \rangle = \langle \theta, \eta\rangle, \quad \theta, \eta \in \cS_{\C}.$$
Using the series expansion \eqref{ser-hol} in the definition of the $\mu_{\beta}$-exponential \eqref{ser-exp}, we get for $\omega \in \cS'$ and $\theta \in \mathcal{U}_{\beta}$
\begin{equation}\label{compare-exp}
e_{\mu_{\beta}}(\omega; \theta) = \sum_{n=0}^{\infty}\frac{1}{n!}\left\langle \sum_{k=0}^{\left\lfloor n/2\right\rfloor}\frac{b_k n!}{2^k (n-2k)!}\tau^{\otimes k} \otimes \omega^{\otimes (n-2k)}, \theta^{\otimes n}\right\rangle.
\end{equation}
}

We define
$$ \mathbb{P}^{\mu_{\beta}}:=\left\{ \left\langle P_{n}^{\mu_{\beta}}(\omega), \theta^{(n)}\right\rangle \middle| \theta^{(n)}\in \cS_{\C}^{\widehat{\otimes} n}, n\in \N_0 \right\}$$
and call it $\mathbb{P}^{\mu_{\beta}}$-system of polynomials on $\cS'_{\C}$.

\rmk{Note that the Appell polynomials are defined by the Taylor series of the $\mu_{\beta}$-exponential \cite{GJRdS15}, i.e.
\begin{equation}\label{mlf-exp-series}
e_{\mu_{\beta}}\left(\omega; \theta \right)=\sum_{n=0}^{\infty}\frac{1}{n!}\left\langle P_{n}^{\mu_{\beta}}(\omega), \theta^{\otimes n}\right\rangle, \quad \omega\in \cS'_{\C}, \theta \in \cU_{0},
\end{equation}
where the kernels $P_{n}^{\mu_{\beta}}(\omega)\in (\cS_{\C}^{\widehat{\otimes} n})'$. The monomial $\left\langle P_{n}^{\mu_{\beta}}(\omega), \theta^{\otimes n}\right\rangle$ may be extended to general kernels $\theta^{(n)}\in \cS_{\C}^{\widehat{\otimes} n}$.}

A comparison with the expansion of the $\mu_{\beta}$-exponential \eqref{compare-exp}, the Appell polynomials with respect to the measure $\mu_{\beta}$ have an explicit formula given by
\begin{equation}\label{Pn-explicit}
P_n^{\mu_{\beta}}(\omega)=\sum_{k=0}^{\left\lfloor n/2\right\rfloor}\frac{b_k n!}{2^k (n-2k)!}\tau^{\otimes k}\otimes \omega^{\otimes (n-2k)}, \omega\in \cS',
\end{equation} 
where the coefficients $b_k$, $k\in \mathbb{N}$ are determined by the recursion formula
$$b_n=-\sum_{k=1}^n\frac{b_{n-k}}{\Gamma(\beta k+1)}, \quad n>0$$
and $b_0=1$.


In the following, we gather various properties of the polynomials $P_n^{\mu_{\beta}}(\cdot)$.
\prop{ \label{appell-prop}\cite{KSWY98} The system of Appell polynomials has the following properties:
\begin{enumerate}
\item[(P1)] For any $x\in \cN'$ and $n \in \mathbb{N}$
\begin{equation}
P_{n}^{\mu_{\beta}}(x)=\sum_{k=0}^{n}\binom{n}{k}x^{\otimes k} \widehat{\otimes}P_{n-k}^{\mu_{\beta}}(0),
\end{equation}
\item[(P2)] For every $x\in \cN'$ and the moment kernels $M_{n}^{\mu_{\beta}}\in \cN'^{\widehat{\otimes}n}$ of $\mu_{\beta}$ defined as follows
\begin{align*}
    \langle M_{n}^{\mu_{\beta}}, \theta^{\otimes n} \rangle &= \langle M_{n}^{\mu_{\beta}}, \theta_1 \widehat{\otimes} \cdots \widehat{\otimes} \theta_n  \rangle \\
    &= \frac{\partial^n}{\partial t_1 \cdots \partial t_n}l_{\mu_{\beta}}(t_1 \theta_1 + \cdots + t_n \theta_n) \vert_{t_1=\cdots=t_n=0}
\end{align*}
we have
\begin{equation}
x^{\otimes n}=\sum_{k=0}^{n}\binom{n}{k}P_{k}^{\mu_{\beta}}(x)\widehat{\otimes}M_{n-k}^{\mu_{\beta}}.
\end{equation}
\item[(P3)] \label{Appell-P3} For all $x,y\in \cN'$ and $M_{n}^{\mu_{\beta}}\in \cN'^{\widehat{\otimes}n}$ as above
\begin{align*}
P_{n}^{\mu_{\beta}}(x+y)
&=\sum_{k+l+m=n}\frac{n!}{k!l!m!}P_{k}^{\mu_{\beta}}(x)\widehat{\otimes}P_{l}^{\mu_{\beta}}(y)\widehat{\otimes} M_{m}^{\mu_{\beta}}\\
&=\sum_{k=0}^{n}\binom{n}{k}P_{k}^{\mu_{\beta}}(x)\widehat{\otimes}y^{\otimes (n-k)}
\end{align*}
\item[(P4)] Further, we observe
$$\mathbb{E_{\mu_{\beta}}}(\langle P_{m}^{\mu_{\beta}}(\cdot), \varphi^{(m)}\rangle)=0 \quad \text{for $m\neq 0$, $\varphi^{(m)}\in \cN_{\C}^{\widehat{\otimes}m}$.}$$
\item[(P5)] For all $p> p_0$ such that the embedding $\iota_{p,p_0}:\cH_p \hookrightarrow \cH_{p_0}$ is a Hilbert-Schmidt and for all $\varepsilon>0$ small enough $\left( \varepsilon\leq \frac{2^{-q_0}}{e \|\iota_{p,p_0}\|_{HS}}\right)$ there exist a constant $C_{p,\varepsilon}>0$ with 
\begin{equation*}
|P_n^{\mu_{\beta}}(\omega)|_{-p}\leq C_{p, \varepsilon} n! \varepsilon^{-n}e^{\varepsilon|\omega|_{-p}}, \quad \omega\in \cH_{-p, \C}.
\end{equation*}
\end{enumerate}}

We now construct a system of generalized functions such that it becomes orthogonal to the $\mathbb{P}^{\mu_{\beta}}$-system above. To this end, we consider the triple
$$\cP(\cS') \subset L^2(\mu_{\beta})\subset \cP'_{\mu_{\beta}}(\cS').$$
The dual pairing between $\cP(\cS')$ and $\cP'_{\mu_{\beta}}(\cS')$ is an extension of the scalar product on $L^2(\mu_{\beta})$, that is, 
$$ \left\langle \left\langle f,\varphi  \right\rangle \right\rangle_{\mu_{\beta}}=(f, \varphi)_{L^2(\mu_{\beta})}, \quad \text{$\varphi \in \cP(\cS')$, $f\in L^2(\mu_{\beta}))$.}$$

For each $\Phi\in \cP'_{\mu_{\beta}}(\cS')$, there is a unique sequence of kernels $\left(\Phi^{(n)}\right)_{n\in \N}\subset \left(\cS_{\C}^{\widehat{\otimes}n}\right)'$, such that
$$\Phi = \sum_{n=0}^{\infty}Q_{n}^{\mu_{\beta}}\left(\Phi^{(n)}\right).$$
Conversely, every series of this type defines a generalized function in $\cP(\cS')$. We have the following system
$$\mathbb{Q}^{\mu_{\beta}}:=\left\{ Q_{n}^{\mu_{\beta}}(\Phi^{(n)})|\Phi^{(n)}\in \left(\cS_{\C}^{\widehat{\otimes}n}\right)', n\in \N_0 \right\}$$
and we call it the $\mathbb{Q}^{\mu_{\beta}}$-system. The pair $(\mathbb{P}^{\mu_{\beta}}, \mathbb{Q}^{\mu_{\beta}})$ is the Appell system and this satisfies the biorthogonality condition:

\thm{\cite{Jah15} \label{ortho-cond-system}
For every $\Phi^{(n)}\in \left(\cS_{\C}^{\widehat{\otimes}n}\right)'$ and $\theta^{(m)} \in \cS_{\C}^{\widehat{\otimes}m}$
$$ \left\langle \! \left\langle \left\langle P_{m}^{\mu_{\beta}}(\cdot), \theta^{(m)}\right\rangle,Q_{n}^{\mu_{\beta}}(\Phi^{(n)})  \right\rangle \! \right\rangle_{\mu_{\beta}}= \delta_{nm}n!\left\langle \Phi^{(n)}, \theta^{(n)} \right\rangle, \quad n,m\in \N_0.$$}

We define a differential operator of order n on $\mathcal{P}(\cS')$ with constant coeffiecient $\Phi^{(n)}\in \left(\cS_{\C}^{\widehat{\otimes}n}\right)'$ by
$$D(\Phi^{(n)})\left\langle \omega^{\otimes m}, \varphi^{(m)}  \right\rangle:= \begin{cases}
\frac{m!}{(m-n)!} \left\langle \omega^{\otimes (m-n)} \widehat{\otimes}\Phi^{(n)} , \varphi^{(m)}\right\rangle & m\geq n, \\
0, & m<n
\end{cases}
$$
for a monomial $\omega \mapsto \left\langle \omega^{\otimes m}, \varphi^{(m)}\right\rangle$ with $\varphi^{(m)}\in \cS_{\C}^{\widehat{\otimes}m}$ and put $Q_{n}^{\mu_{\beta}}(\Phi^{(n)})=D(\Phi^{(n)})^{\ast}1$.

\lem{\cite{KSWY98} $D(\Phi^{(n)})$ is a continous linear operator from $\cP(\cS')$ to $\cP(\cS')$.}

\rmk{For $D(\Phi^{(1)})\in \cS'$, we have the usuaal G${\hat{a}}$teaux derivative as e.g. in white noise analysis \cite{HKPS93}
\begin{equation}
D(\Phi^{(1)})\varphi= D_{\Phi^{(1)}}\varphi := \left. \frac{d}{dt}\varphi(\cdot + t\Phi^{(1)})\right|_{t=0}
\end{equation}
for $\varphi\in \cP(\cS)$ and we have $D\left((\Phi^{(1)})^{\otimes n}\right)=\left(D_{\Phi^{(1)}}\right)^n$ thus $D\left((\Phi^{(1)})^{\otimes n}\right)$ is in fact a differential operator of order $n$.}

\rmk{\cite{KSWY98}
For $\Phi^{(n)}\in \left(\cS_{\C}^{\widehat{\otimes}n}\right)'$, $\varphi^{(m)}\in \cS_{\C}^{\widehat{\otimes}m}$, we have
\begin{equation}
D(\Phi^{(n)})\left\langle P_m^{\mu_{\beta}}(x), \varphi^{(m)}  \right\rangle:= \begin{cases}
\frac{m!}{(m-n)!} \left\langle P_{m-n}^{\mu_{\beta}}(x) \widehat{\otimes}\Phi^{(n)} , \varphi^{(m)}\right\rangle & m\geq n, \\
0, & m<n.
\end{cases}
\end{equation}}

\noindent With the help of the Appell systems, a test function and a distribution space in non-Gaussian analysis can be constructed, the details of this construction can be found in \cite{KSWY98}, \cite{GJRdS15}, \cite{GJ16} and references therein. In between the many choices of triples which can be constructed, we choose the Kondratiev triple 
\[
(\cS)_{\mu_{\beta}}^{1}\subset(H_{p})_{q,\mu_{\beta}}^{1}\subset L^{2}(\mu_{\beta})\subset(H_{-p})_{-q,\mu_{\beta}}^{-1}\subset(\cS)_{\mu_{\beta}}^{-1}.
\]
We have that the test function space $(\cS)_{\mu_{\beta}}^1$ is the projective limit of the family of spaces $(\cH_p)_{q, \mu_{\beta}}^{1}$, that is,
$$(\cS)_{\mu_{\beta}}^1 = \displaystyle{\bigcap_{p,q\in \N_0}}(\cH_p)_{q, \mu_{\beta}}^{1}$$
and $(\cS)_{\mu_{\beta}}^1$ is equipped with the coarser topology such that the embeddings
$$(\cS)_{\mu_{\beta}}^1 \hookrightarrow (\cH_p)_{q, \mu_{\beta}}^{1}$$
are continuous. The test function space $(\cS)_{\mu_{\beta}}^1$ is a nuclear space which is continuously embedded in $L^2(\mu_{\beta})$. 
By the general duality theory, the space of generalized functions $(\cS)_{\mu_{\beta}}^{-1}$ is the dual space of $(\cS)_{\mu_{\beta}}^{1}$ with respect to $L^2(\mu_{\beta})$ and is defined as the inductive limit of the family $(\cH_{-p})_{-q, \mu_{\beta}}^{-1}$, that is,
$$(\cS)_{\mu_{\beta}}^{-1}=\displaystyle{\bigcup_{p,q\in \N_0}}(\cH_{-p})_{-q, \mu_{\beta}}^{-1}.$$
The space $(\cS)_{\mu_{\beta}}^{-1}$ is equipped with the finest topology such that the embeddings
$$(\cH_{-p})_{-q, \mu_{\beta}}^{-1} \hookrightarrow (\cS)_{\mu_{\beta}}^{-1}.$$
The space $(H_{p})_{q,\mu_{\beta}}^{1}$ is defined as the completion of the $\mathcal{P}\big(S'\big)$ (the space of smooth polynomials on $S'$) w.r.t.\ the norm $\|\cdot\|_{p,q,\mu_{\beta}}$ given by 
\[
\|\varphi\|_{p,q,\mu_{\beta}}^{2}:=\sum_{n=0}^{\infty}(n!)^{2}2^{nq}|\varphi^{(n)}|_{p}^{2},\quad p,q\in\mathbb{N}_{0},\;\varphi\in\mathcal{P}\big(S'\big).
\]
The dual space $(H_{-p})_{-q,\mu_{\beta}}^{-1}$ is a subset of $\mathcal{P}'\big(S'\big)$ such that if $\Phi\in(H_{-p})_{-q,\mu_{\beta}}^{-1}$, then 
\[
\|\Phi\|_{-p,-q,\mu_{\beta}}^{2}:=\sum_{n=0}^{\infty}2^{-nq}|\Phi^{(n)}|_{-p}^{2}<\infty,\quad p,q\in\mathbb{N}_{0}.
\]
The dual pairing between $\big(S'\big)_{\mu_{\beta}}^{-1}$ and $\big(S\big)_{\mu_{\beta}}^{1}$,
denoted by $\langle\!\langle\cdot,\cdot\rangle\!\rangle_{\mu_{\beta}}$
is a bilinear extension of scalar product in $L^{2}(\mu_{\beta})$.

The action of a distribution 
\begin{equation}\Phi = \sum_{n=0}^{\infty}Q_{n}^{\mu_{\beta}}(\Phi^{(n)}) \in  (\cS)_{\mu_{\beta}}^{-1}
\end{equation}
on a test function
\begin{equation}
    \label{tp}\varphi = \sum_{n=0}^{\infty}\left\langle P_{n}^{\mu_{\beta}}(\omega), \varphi^{(n)} \right\rangle \in (\cS)_{\mu_{\beta}}^{1}
\end{equation}
using the biorthogonality property in Theorem \ref{ortho-cond-system} is given by 
$$\left\langle \left\langle \Phi, \theta \right\rangle \right\rangle=\sum_{n=0}^{\infty}n! \left\langle \Phi^{(n)}, \theta^{(n)}\right\rangle.$$

The set of $\mu_{\beta}$-exponentials 
\[
\left\{ e_{\mu_{\beta}}(\cdot,\varphi):=\frac{e^{\langle\cdot,\varphi\rangle}}{\mathbb{E}\big(e^{\langle\cdot,\varphi\rangle}\big)},\;\varphi\in S_{\mathbb{C}},\;|\varphi|_{p}<2^{-q}\right\} 
\]
forms a total set in $(H_{p})_{q,\mu_{\beta}}^{1}$ and for any $\varphi\in S_{\mathbb{C}}$
such that $|\varphi|_{p}<2^{-q}$ we have $\|e_{\mu_{\beta}}(\cdot,\varphi)\|_{p,q,\mu_{\beta}}<\infty$.

Let us introduce an integral transform, the $S_{\mu_{\beta}}$-transform,
which is used to characterize the spaces $(S)_{\mu_{\beta}}^{1}$
and $(S)_{\mu_{\beta}}^{-1}$. For any $\Phi\in(S)_{\mu_{\beta}}^{-1}$
and $\varphi\in U\subset S_{\mathbb{C}}$, where $U$
is a suitable neighborhood of zero, we define 
\[
S_{\mu_{\beta}}\Phi(\varphi):=\frac{\langle\!\langle\Phi,e^{\langle\cdot,\varphi\rangle}\rangle\!\rangle_{\mu_{\beta}}}{\mathbb{E}\big(e^{\langle\cdot,\varphi\rangle}\big)}=\frac{1}{E_{\beta}(\frac{1}{2}\langle\varphi,\varphi\rangle)}\langle\!\langle\Phi,e^{\langle\cdot,\varphi\rangle}\rangle\!\rangle_{\mu_{\beta}}.
\]
The characterization theorem for the space $(S)_{\mu_{\beta}}^{-1}$
via the $S_{\mu_{\beta}}$-transform is done using the spaces of holomorphic
functions on $S_{\mathbb{C}}$. Note that the space $\mathrm{Hol}_{0}\big(S_{\mathbb{C}}\big)$
is given as the inductive limit of a family of normed spaces, see
\cite{KSWY98} for the details and the proof of the following characterization
theorem. 
\thm{[{cf.\ \cite[Theorem~8.34]{KSWY98}}]
 \label{theorem:Charact_distributions}The $S_{\mu_{\beta}}$-transform
is a topological isomorphism from $(S)_{\mu_{\beta}}^{-1}$ to $\mathrm{Hol}_{0}\big(S_{d,\mathbb{C}}\big)$.
}

As a corollary from the characterization theorem the following integration result can be deduced. For details and proofs we refer to \cite{GJRdS15}.
\thm{\label{it}
Let $(T,\mathcal{B},\nu)$ be a measure space and $\Phi_t\in (S)_{\mu_{\beta}}^{-1}$ for all $t\in T$. Let $\mathcal{U}\subset S_{\mathbb{C}}$ be an appropriate neighbourhood of zero and $0<C<\infty$, such that
\begin{enumerate}
\item $S_{\mu_{\beta}}\Phi_{\cdot}(\xi):T\to \mathbb{C}$ is measurable for all $\xi \in \mathcal{U}$.
\item $\int_T\left| S_{\mu_{\beta}}\Phi_t(\xi)\right|\, d\nu(t)\le C$ for all $\xi\in \mathcal{U}$.
\end{enumerate}
Then, there exists $\Phi\in (S)_{\mu_{\beta}}^{-1}$ such that for all $\xi\in \mathcal{U}$
\[ S_{\mu_{\beta}}\Psi(\xi)=\int_TS_{\mu_{\beta}}\Phi_t(\xi)\, d\nu(t). \]
We denote $\Psi$ by $\int_T\Phi_t\, d\nu(t)$ and call it the weak integral of $\Phi$.
}
\\

In the following we will use the $T_{\mu_{\beta}}$-transform which is defined as follows.
\lem{
Let $\Phi \in (S)_{\mu_{\beta}}^{-1}$ and $p,q\in \mathbb{N}$ such that $\Phi \in  (H_{-p})_{-q,\mu_{\beta}}^{-1}$. Then, the $T_{\mu_{\beta}}$-transform given by
\[ T_{\mu_{\beta}}\Phi(\varphi)=\langle\!\langle \Phi, \exp\left( i\langle \cdot,\varphi\rangle\right) \rangle\!\rangle_{\mu_{\beta}} \]
is well-defined for $\varphi\in U_{p,q}$ and we have
\[  T_{\mu_{\beta}}\Phi(\varphi)=E_{\beta}\left(-\frac{1}{2}\langle \varphi,\varphi\rangle \right)S_{\mu_{\beta}}\Phi(i\varphi).\]
In particular, $T_{\mu_{\beta}}\Phi\in \mathrm{Hol}_{0}\big(S_{\mathbb{C}}\big)$ if and only if $S_{\mu_{\beta}}\in \mathrm{Hol}_{0}\big(S_{\mathbb{C}}\big)$. Moreover, Theorem \ref{it} also holds if the $S_{\mu_{\beta}}$-transform is replaced by the $T_{\mu_{\beta}}$-transform.}\\

\section{Symbol of an operator}
\defn{\label{sym-defn} For $\Xi \in \cL\left((\cH_p)_q^1, (\cS)^{-1} \right)$, we define the \emph{operator symbol} of $\Xi$ by
\begin{equation}
    \widehat{\Xi}(\xi, \eta)=\langle \! \langle  \Xi e_{\mu_{\beta}}(\cdot, \xi), e_{\mu_{\beta}}(\cdot, \eta) \rangle \! \rangle, \quad \xi, \eta \in \cS_{\C}(\R),
\end{equation}
with $|\xi|_p, |\eta|_p< \frac{C}{2^{q}}$ for $C, q>0$.}

We denote the set $E_{p,q}=\{\xi, \eta :  |\xi|_p, |\eta|_p< \frac{C}{2^{q}}$ for $C, q>0\}$.

\rmk{The set $A_{p,q}=\{e_{\mu_{\beta}}(\cdot, \xi) \vert |\xi|_p<\frac{C}{2^q}\}$ is total in $(\cH_p)_q^1$.}

The $S_{\mu_{\beta}}$-transform of $\Phi\in (\cS)^{-1}$ is defined by 
\begin{equation*}
    S_{\mu_{\beta}}\Phi(\xi)= \langle\! \langle \Phi, e_{\mu_{\beta}}(\cdot; \xi)\rangle\!\rangle, \quad \xi \in \cS_{\C}(\R).
\end{equation*}
We have the following:
\rmk{\label{rem-sym}For $\Xi \in \cL\left((\cH_p)_q^1, (\cS)^{-1} \right)$, it holds that
\begin{equation}
 \widehat{\Xi}(\xi, \eta) = S_{\mu_\beta}\left( \Xi e_{\mu_{\beta}}(\cdot, \xi) \right)(\eta) = S_{\mu_\beta}\left( \Xi^{*} e_{\mu_{\beta}}(\cdot, \eta) \right)(\xi), \quad \xi, \eta \in \cS_{\C}(\R).
\end{equation}}

We now discuss analytic properties of the symbol of an operator $\Xi \in \cL\left((\cH_p)_q^1, (\cS)^{-1} \right)$.

\cor{ \label{entireness-sym}If $\Xi \in \cL\left((\cH_p)_q^1, (\cS)^{-1} \right)$, then for any $\xi_1, \xi_2, \eta_1, \eta_2 \in E_{p,q}$, we have
\begin{equation*}
    z, \omega \mapsto \widehat{\Xi}(\xi_1+z\xi_2, \eta_1+\omega\eta_2), \quad z, \omega \in \C, 
\end{equation*}
is an entire holomorphic function on $\C\times\C$.
}

\pf{Let $\Xi \in \cL\left((\cH_p)_q^1, (\cS)^{-1} \right)$ and $\xi_1, \xi_2, \eta_1, \eta_2 \in E_{p,q}$. Further, for $e_{\mu_{\beta}}(\cdot, \xi_1+z\xi_2)$ given by the series expansion
\begin{equation*}
    e_{\mu_{\beta}}(\cdot, \xi_1+z\xi_2)=\sum_{n=0}^{\infty}\frac{1}{n!}\langle P_n^{\mu_{\beta}}, (\xi_1+z\xi_2)^{\otimes n} \rangle,
\end{equation*}
$\Xi e_{\mu_{\beta}}(\cdot, \xi_1+z\xi_2)\in L^2(\mu_{\beta})\subset (\cS)^{-1}$. \\
Then by the characterization theorem \ref{theorem:Charact_distributions}, $S_{\mu_\beta}\left( \Xi e_{\mu_{\beta}}(\cdot, \xi_1+z\xi_2) \right)(\eta_1+\omega\eta_2)$ is entire holomorphic in $\omega$. Similar argument follows for $S_{\mu_\beta}\left( \Xi^{*} e_{\mu_{\beta}}(\cdot, \eta_1+\omega\eta_2) \right)$, this means that $$ S_{\mu_\beta}\left( \Xi^{*} e_{\mu_{\beta}}(\cdot, \eta_1+\omega\eta_2) \right)(\xi_1+z\xi_2)$$ is entire holomorphic in $z$.\\
Using the fact that 
\begin{align*}
\widehat{\Xi}(\xi_1+z\xi_2, \eta_1+\omega\eta_2) &= S_{\mu_\beta}\left( \Xi e_{\mu_{\beta}}(\cdot, \xi_1+z\xi_2) \right)(\eta_1+\omega\eta_2) \\
&= S_{\mu_\beta}\left( \Xi^{*} e_{\mu_{\beta}}(\cdot, \eta_1+\omega\eta_2) \right)(\xi_1+z\xi_2), \quad \xi, \eta \in \cS_{\C}(\R),
\end{align*}
by Remark \ref{rem-sym}. Accordingly, $\widehat{\Xi}(\xi_1+z\xi_2, \eta_1+\omega\eta_2)$ is entire holomorphic in two variables $z,\omega$. \qed}

\thm{\label{char-thm-op-sym}
For $\xi, \eta \in \cS_{\C}(\R)$, we have $\Xi \in \cL\left((\cH_p)_q^1, (\cS)^{-1} \right)$ if and only if $\widehat{\Xi}(\xi, \eta)\in Hol_{\C}(E_{p,q}\times E_{p,q})$.}

\pf{The proof for necessity part follows from Corollary \ref{entireness-sym}. 

Now, let $\Xi \in \cL\left((\cH_p)_q^1, (\cS)^{-1} \right)$, then for all $\xi, \eta\in (\cH_{q})_p^{1}$, $\Xi e_{\mu_{\beta}}(\cdot, \xi) \in (\cS)^{-1}$ and $\Xi e_{\mu_{\beta}}(\cdot, \eta) \in (\cS)^{-1}$. Hence, $
 S_{\mu_\beta}\left( \Xi e_{\mu_{\beta}}(\cdot, \xi) \right)(\eta)$ and $S_{\mu_\beta}\left( \Xi e_{\mu_{\beta}}(\cdot, \eta) \right)(\xi) \in Hol_{\C}(E_{p,q})$ by characterization theorem \ref{theorem:Charact_distributions} in $\eta$, and $\xi$ respectively. Accordingly, $\widehat{\Xi}(\xi, \eta)\in Hol_{\C}(E_{p,q}\times E_{p,q})$.\qed}

 \cor{\label{char-thm-op-sym-ext}
For $\xi, \eta \in \cS_{\C}(\R)$, we have $\Xi \in \cL\left((\cH_p)_q^1, (\cS)^{-1} \right)$ if and only if $\widehat{\Xi}(\xi, \eta)\in Hol_{\C}(E_{p,q}\times \cS(\R))$.}

\ex{\label{char-diff-ex}
Consider $D_{\psi}$ for $\psi \in L^2(\R)$. Taking the operator symbol of $D_{\psi}$ we have
\begin{align*}
\widehat{D}_{\psi}(\xi, \eta) &= \langle \! \langle D_{\psi}e_{\mu_{\beta}}(\cdot, \xi), e_{\mu_{\beta}}(\cdot, \eta)  \rangle \! \rangle \\
&= \langle \psi, \xi \rangle \langle \! \langle e_{\mu_{\beta}}(\cdot, \xi), e_{\mu_{\beta}}(\cdot, \eta)  \rangle \! \rangle.
\end{align*}
Now, observe that
\begin{align*}
\widehat{D}_{\psi}(\xi_1+z\xi_{2}, \eta_1 + \omega\eta_2) &= \langle \psi, \xi_1+z\xi_2 \rangle \langle \! \langle e_{\mu_{\beta}}(\cdot, \xi_1+z\xi_2), e_{\mu_{\beta}}(\cdot, \eta_1+\omega\eta_2)  \rangle \! \rangle \\
&= \langle \psi, \xi_1+z\xi_2 \rangle 
 S_{\mu_{\beta}}\left(e_{\mu_{\beta}}(\cdot, \xi_1+z\xi_2) \right) \left(e_{\mu_{\beta}}(\cdot, \eta_1+\omega\eta_2) \right) \\
&= \langle \psi, \xi_1+z\xi_2 \rangle 
 S_{\mu_{\beta}}\left(e_{\mu_{\beta}}(\cdot, \eta_1+\omega\eta_2) \right)\left(e_{\mu_{\beta}}(\cdot, \xi_1+z\xi_2) \right).
\end{align*}
Note that $\langle \psi, \xi_1+z\xi_2 \rangle$ is a polynomial in $z$ hence entire. Further, since $e_{\mu_{\beta}}(\cdot, \xi)\in (\cS)^{-1}$ for all $\xi\in (\cH_{q})_p^{1}$, we have that
$S_{\mu_{\beta}}\left(e_{\mu_{\beta}}(\cdot, \xi_1+z\xi_2) \right) \left(e_{\mu_{\beta}}(\cdot, \eta_1+\omega\eta_2) \right) \in Hol_{\C}(E_{p,q})$ by characterization theorem \ref{theorem:Charact_distributions}. Similar argument follows for $S_{\mu_{\beta}}\left(e_{\mu_{\beta}}(\cdot, \eta_1+\omega\eta_2) \right)\left(e_{\mu_{\beta}}(\cdot, \xi_1+z\xi_2) \right)$. Consequently, $D_{\psi} \in \cL\left((\cH_p)_q^1, (\cS)^{-1} \right)$ by Theorem \ref{char-thm-op-sym}.
}

\section{Some examples}
\subsection{Translation Operator}
\label{rem-trans-vs-multip} Note that for $\eta \in \cS'_{\C}(\R)$, independent of the measure, we have the following observations:
\begin{align*}
    \int_{\cS'}e^{i\langle \omega, \xi\rangle} \frac{e^{\langle \omega, \eta\rangle}}{\mathbb{E}(e^{\langle \cdot, \eta\rangle})}d\mu(\omega) &= \frac{1}{\mathbb{E}(e^{\langle \cdot, \eta\rangle})} \int_{\cS'}{e^{i\langle \omega, \xi-i\eta\rangle}}d\mu(\omega)\\
    &=\frac{C(\xi-i \eta)}{\mathbb{E}(e^{\langle \cdot, \eta\rangle})}, 
\end{align*}
on the other hand,
\begin{align*}
    \int_{\cS'} e^{i \langle \omega+\eta, \xi\rangle} d\mu(\omega) &= C(\xi) \cdot e^{i\langle \eta, \xi \rangle}
\end{align*}
with $C$ as the characteristic function.
Accordingly, both only coincide whenever:
$$\frac{C(\xi)}{C(\xi-i\eta)} = \frac{\mathbb{E}(e^{\langle \cdot, \eta\rangle})}{e^{i\langle \eta, \xi\rangle}}.$$

We define the translation operator (or shift operator) as follows:

\begin{defn}
Let $y \in \cS'$ and $\varphi \in (\cS)$ we define the translation operator $\tau_y$ by $\tau_y \varphi(\omega) = \varphi(\omega + y).$
\end{defn}
See e.g.~\cite{W95, V10}

Now, we consider identity (P3) of the polynomial $P_{n}^{\mu_{\beta}}(\cdot)$ from Proposition \ref{appell-prop}:
\begin{equation}\label{shift-mlf}
     P_{n}^{\mu_{\beta}}(x+y)=\sum_{k=0}^{n}\binom{n}{k}P_{k}^{\mu_{\beta}}(x)\widehat{\otimes}y^{\otimes (n-k)}.  
\end{equation}
The next proposition gives expansion of the operator:
\prop{
The translation operator $\tau_y$ for $y \in \cS'$ can be represented as
\begin{equation*}
    \tau_y \varphi(\omega)= \varphi(\omega+y)=
    \sum_{k=0}^{\infty} \sum_{n=0}^{\infty} \binom{n+k}{k} \left\langle P_k^{\mu_{\beta}}(\omega) , y^{\otimes n}\widehat{\otimes}_{n} \varphi^{(n+k)}\right\rangle, \quad \varphi \in (\cS)^1.
\end{equation*}
Furthermore, for any $p\geq 0$, $q>0$ with $|y|_{-(p+q)}<\infty$, it holds that 
\begin{equation}
    \|\tau_{y}\varphi\|_p \leq \frac{\| \varphi\|_{p+q}}{(1-2^{2q})^{1/2}}{\rm{exp}}\left( \frac{|y|^2_{-(p+q)}}{2(1-2^{2q})}\right).
\end{equation}
}

\pf{For $y \in \cS', \varphi \in (\cS)^1$ with kernels $\varphi^{(n)}$, $n\in \mathbb{N}$, i.e.,
$$\varphi(\omega)=\sum_{n=0}^{\infty}\langle P_n^{\mu_{\beta}}(\omega), \varphi^{(n)}\rangle.$$
and in view of identity \eqref{shift-mlf}, we obtain
\begin{align*}
    \tau_y \varphi(\omega)&= \varphi(\omega+y) \\
    &= \sum_{n=0}^{\infty}\langle P_n^{\mu_{\beta}}(\omega+y), \varphi^{(n)}\rangle \\
    &= \sum_{n=0}^{\infty}\left\langle \sum_{k=0}^n \binom{n}{k}P_k^{\mu_{\beta}}(\omega) \widehat{\otimes} y^{\otimes (n-k)}, \varphi^{(n)}\right\rangle, \tag*{by Proposition 15} \\
    &= \sum_{n=0}^{\infty} \sum_{k=0}^n \binom{n}{k} \left\langle P_k^{\mu_{\beta}}(\omega) , y^{\otimes (n-k)}\widehat{\otimes}_{n-k} \varphi^{(n)}\right\rangle \\
    &= \sum_{k=0}^{\infty} \sum_{n=0}^{\infty} \binom{n+k}{k} \left\langle P_k^{\mu_{\beta}}(\omega) , y^{\otimes n}\widehat{\otimes}_{n} \varphi^{(n+k)}\right\rangle
\end{align*}
Using the inequality
$$ |y^{\otimes n}\widehat{\otimes}_{n} \varphi^{(n+k)} |_p \leq 2^{qk} |y|^n_{-(p+q)} |\varphi^{(n-k)}|_{p+q}$$
we have
\begin{align*}
    \|\tau_y \varphi\|_p^2 &=\sum_{k=0}^{\infty} k! \left|\sum_{n=0}^{\infty}\binom{n+k}{k}  y^{\otimes n}\widehat{\otimes}_{n} \varphi^{(n+k)}\right|_p^2 \\
    &\leq \sum_{k=0}^{\infty} k! \left( \sum_{n=0}^{\infty} \frac{(n+k)!}{n!k!} 2^{qk} |y|^n_{-(p+q)} |\varphi^{(n-k)}|_{p+q} \right)^2 \\
    &\leq \sum_{k=0}^{\infty} k! \left( \sum_{n=0}^{\infty} (n+k)!  |\varphi^{(n-k)}|_{p+q}^2\right) \left( \sum_{n=0}^{\infty} \frac{(n+k)!}{(n!k!)^2} 2^{2qk} |y|^{2n}_{-(p+q)}\right) \\
    &\leq \|\varphi\|_{p+q}^2 \sum_{n=0}^{\infty}\sum_{k=0}^{\infty} \frac{(n+k)!}{n!n!k!} 2^{2qk} |y|^{2n}_{-(p+q)}.
\end{align*}
Since 
$$ \sum_{k=0}^{\infty} \frac{(n+k)!}{n!k!} 2^{2qk}= (1-2^{2q})^{-(n+1)},$$
we conclude that 
\begin{align*}
\|\tau_y \varphi\|_p^2 &\leq \|\varphi\|_{p+q}^2 \sum_{n=0}^{\infty} \frac{1}{n!}|y|^{2n}_{-(p+q)}  (1-2^{2q})^{-(n+1)} \\
&= \|\varphi\|_{p+q}^2 (1-2^{2q})^{-1} {\rm{exp}} \left( |y|^2_{-(p+q)}(1-2^{2q})^{-1}\right).
\end{align*}
Consequently, we have 
\begin{align*}
        \|\tau_{y}\varphi\|_p \leq \frac{\| \varphi\|_{p+q}}{(1-2^{2q})^{1/2}}{\rm{exp}}\left( \frac{|y|^2_{-(p+q)}}{2(1-2^{2q})}\right). 
\end{align*} \qed
}

\rmk{\label{trans-exp-mlf}
With $\mu_{\beta}$-exponential given by
\begin{equation*}
e_{\mu_{\beta}}\left(\omega, \theta \right)= \frac{e^{\left(\omega, \theta \right)}}{E_\beta\left(\frac{1}{2}\langle \theta, \theta\rangle\right)},
\end{equation*}
the action of translation operator to the $e_{\mu_{\beta}}$-exponential is as follows
\begin{align*}
    \tau_{\eta}e_{\mu_{\beta}}\left(\omega, \theta \right) &= \tau_{\eta}\left(\frac{e^{\left(\omega, \theta \right)}}{E_\beta\left(\frac{1}{2}\langle \theta, \theta\rangle\right)} \right) \\
    &=\frac{e^{\left(\omega+\eta, \theta \right)}}{E_\beta\left(\frac{1}{2}\langle \theta, \theta\rangle\right)} \\
    &=\frac{e^{\left(\omega, \theta \right)}e^{\left(\eta, \theta \right)}}{E_\beta\left(\frac{1}{2}\langle \theta, \theta\rangle\right)} \\
    &=e^{\left(\eta, \theta \right)} e_{\mu_{\beta}}\left(\omega, \theta \right), \quad \eta\in \cS', \theta\in \cS.
\end{align*}
}

If we use the definition of the corresponding differential operator on the polynomials, test functions and exponential vectors, we obtain: 
\begin{eqnarray*}
    \langle\! \langle  \tau_{\eta} \exp(i \langle \cdot, \xi), \exp(i \langle \cdot, \varphi) \rangle \! \rangle &=& \langle\! \langle \exp(i \langle \cdot + \eta, \xi),\exp(i \langle \cdot, \varphi) \rangle \! \rangle\\
    &=& \exp(i \langle \eta, \xi \rangle )T_{\mu_{\beta}}(\exp(i \langle \cdot, \xi)) (\varphi) 
\end{eqnarray*}
on the other hand: 
\begin{eqnarray*}
    \langle\! \langle \exp(D_{\eta})  \exp(i \langle \cdot, \xi), \exp(i \langle \cdot, \varphi) \rangle \! \rangle
    &=&  \left\langle\! \left\langle  \sum_{k=0}^{\infty} \frac{1}{k!} D^k_{\eta} \exp(i \langle \cdot, \xi), \exp(i \langle \cdot, \varphi) \right\rangle \! \right\rangle \\
    &=& \left\langle\! \left\langle \sum_{k=0}^{\infty} \frac{(i)^k}{k!} \langle \eta, \xi\rangle^k  \exp(i \langle \cdot, \xi), \exp(i \langle \cdot, \varphi) \right\rangle \! \right\rangle \\
    &=& \exp(i \langle \eta, \xi \rangle )T_{\mu_{\beta}}(\exp(i \langle \cdot, \xi)) (\varphi)
\end{eqnarray*}
For polynomials we have a similar derivation. We hence have the following theorem.

\thm{For any test function $\Phi \in (S)^1$, we have that $\tau_{\eta}= \exp(D_{\eta}),\quad \eta \in \cS_{\mathbb{C}}$. Moreover, this equality can be generalized to infinitely often differentiable functions with kernels in the Schwartz test functions. }

\ex{
We define the Mehler's formula for the Ornstein-Uhlenbeck semigroup in the Gaussian case. Indeed, we have for all $\Phi \in \cS$,
\begin{equation}\label{Mehlers-form}
    P_{t}\Phi(y)=\int_{\cS(\R)}\Phi\left( e^{-t} y+ \sqrt{1-e^{-2t}}\omega\right) d\mu(\omega), \quad  y \in \cS'(\R), \,t\geq 0.
\end{equation}
Equivalently, we have 
\begin{equation}\label{Mehlers-form1}
    P_{t}\Phi(y)=\langle \!\langle \tau_{e^{-t} y}\sigma_{\sqrt{1-e^{-2t}}} \Phi, 1\!\!1 \rangle\! \rangle.
\end{equation}
In the case of an exponential function, for $\xi \in \cS(\R)$, we obtain:
\begin{eqnarray*}
P_t e^{i \langle y, \xi \rangle} &=&  
\langle \!\langle \tau_{e^{-t} y}\sigma_{\sqrt{1-e^{-2t}}} e^{i \langle \omega, \xi \rangle}, 1\!\!1 \rangle\! \rangle\\
&=& \langle \!\langle  e^{i \langle e^{-t}y + \sqrt{1-e^{-2t}} \omega, \xi \rangle}, 1\!\!1 \rangle\! \rangle\\
&=& e^{i e^{-t} \langle y, \xi \rangle}\langle \!\langle  e^{i \langle   \omega ,\sqrt{1-e^{-2t}} \xi \rangle}, 1\!\!1 \rangle\! \rangle\\
&=& e^{i e^{-t} \langle y, \xi \rangle} E_{\beta}\left( -\frac{1}{2} (1-e^{-2t}) \langle \xi, \xi\rangle\right).
\end{eqnarray*}
Furthermore, using the integral representation for the Mittag-Leffler function:
\begin{equation}
    E_{\beta}(-z)= \int_{0}^{\infty}e^{-m z}d\nu_{\beta}(m),
\end{equation}
where $\nu_{\beta}$ is a probability measure on $[0, \infty)$ with density $M_{\beta}$ with respect to the Lebesgue measure, that is, $d\nu_{\beta}(m)= M_{\beta}(m)d m$ yields the following:
\begin{eqnarray*}
P_t e^{i \langle y, \xi \rangle} 
&=& e^{i e^{-t} \langle y, \xi \rangle} E_{\beta}\left( -\frac{1}{2} (1-e^{-2t}) \langle \xi, \xi\rangle\right) \\
&=& e^{i \langle y, e^{-t} \xi \rangle} \int_{0}^{\infty} e^{-m \left( \frac{1}{2}(1-e^{-2t})\|\xi\|^2\right)} M_{\beta}(m)d m \\
&=& \int_{0}^{\infty} e^{-\frac{1}{2} \left( m(1-e^{-2t})\|\xi\|^2+i \langle y, e^{-t}\xi \rangle\right)} M_{\beta}(m)d m \\ 
&=& \int_{0}^{\infty} \mathbb{E} \left( \tau_{e^{-t} y}\sigma_{\sqrt{1-e^{-2t}}} e^{i \langle \sqrt{m}\omega, \xi \rangle} \right)M_{\beta}(m)d m.
\end{eqnarray*}
It is quite easy to prove that $(P_t)_{t\geq 0}$ is not having the semi-group property. 

\subsection{Integral kernel operators}
In this section, we begin with recalling a differential operator $\partial_t$ which plays a fundamental role in the white noise calculus. For $t \in T$ and $f\in \cS_{\C}^{\widehat{\otimes} (n+1)}$, we define $\delta_t \widehat{\otimes}_1 f \in \cS_{\C}^{\widehat{\otimes} n}$ by 
\begin{equation*}
    \delta_{t}\widehat{\otimes}_1 f(t_1, \dots, t_n)=f(t,t_1, \dots , t_n), \quad t_1, \dots, t_n \in T.
\end{equation*}

\thm{\label{diff-testfunc}
Assume that $\varphi \in (\cS)^1$ is given as in \eqref{tp}. Then the series
\begin{equation}\label{diff-ope}
D_{y}\varphi(\omega)=\sum_{n=0}^{\infty}n\left\langle P_{n-1}^{\mu_{\beta}}(\omega), y \widehat{\otimes}_1\varphi^{(n)} \right\rangle, \quad \omega,y\in \cS'
\end{equation}
converges absolutely. Moreover, $D_{y}\varphi \in (\cS^1)$ and, for any $p\geq 0$ and $q>0$ we have
\begin{equation}
\|D_{y}\varphi\|_{p,q}\leq \left(\frac{2^{-2q}}{-2qe\log(2)}\right)^{1/2}|y|_{-(p+q)}\|\varphi\|_{p+q}, \quad \varphi\in(\cS^1).
\end{equation} 
In particular, $D_{y}\in \cL((\cS^1), (\cS^1))$.}

If we take $y=\delta_t$, the Dirac delta function at $t$, then
\begin{enumerate}
\item $\partial_t=D_{\delta_t}$ is called the white noise differential operator, or the Hida differential operator, or the annihilation operator
\item $\partial_t^*=D_{\delta_t}^*$ is called the creation operator.
\end{enumerate}
Having introduced the differential operator $\partial_t$ above, we now present a general theory of operators which are expressed as an integral of $\partial_t$ and $\partial_t^*$.

We introduce an operator in $\cL((\cS)^{1}, (\cS)^{-1})$ which is expressed in a formal integral:
\begin{equation}
    \int_{T^{l+m}} \kappa(s_1, \dots, s_{l}, t_1, \dots, t_m)\partial_{s_1}^*\cdots \partial_{s_l}^*\partial_{t_1}\cdots \partial_{t_m}ds_1\cdots ds_l dt_1\cdots dt_m, 
\end{equation}
where $\kappa \in (\cS_{\C}^{\otimes (l+m)})^*$.
\lem{\label{int-ker} For $\varphi \in (\cS)^1$ and $\psi \in (\cS)^{-1}$, we put
    \begin{equation}
    \eta_{\varphi, \psi} (s_1,\dots, s_l, t_1, \dots, t_m) = \langle \! \langle \psi, \partial_{s_1}^*\cdots \partial_{s_l}^* \partial_{t_1} \cdots\partial_{t_m} \varphi \rangle \! \rangle.
   \end{equation}
    Then for any $p>0$ we have 
    $$|\eta_{\varphi, \psi}|_p \leq 2^{-p} \sqrt{l^l m^m} \left(\frac{2^{-p}}{-2pe\log(2)}\right)^{\frac{l+m}{2}} \|\varphi\|_{p,q} \|\psi\|_{p,q}.$$
    In particular $\eta_{\varphi, \psi} \in \cS_{\C}^{\otimes (l+m)}$.}

\thm{\label{int-ker-op}
For any $\kappa\in (\cS_{\C}^{\otimes (l+m)})^{-1}$ there exists a continuous linear operator $\Xi_{l,m}(\kappa)\in \cL((\cS),(\cS)^{-1})$ such that
\begin{equation}
\left\langle\! \left\langle \Xi_{l,m}(\kappa)\varphi, \psi  \right\rangle \!\right\rangle = \langle \kappa, \eta_{\varphi, \psi}\rangle, \quad \varphi, \psi \in (\cS),
\end{equation}
where
 \begin{equation*}
    \eta_{\varphi, \psi} (s_1,\dots, s_l, t_1, \dots, t_m) = \langle \! \langle \partial_{s_1}^*\cdots \partial_{s_l}^* \partial_{t_1} \cdots\partial_{t_m} \varphi, \psi \rangle \! \rangle.
   \end{equation*}
Moreover, for any $p>0$ with $|\kappa|_{-p}<\infty$ it holds that
\begin{equation}\label{intker}
\|\Xi_{l,m}(\kappa)\varphi\|_{-p} \leq 2^{-p}(l^l m^m)^{1/2}\left( \frac{2^{-p}}{-2pe\log{2}}\right)^{(l+m)/2}|\kappa|_{-p}\|\varphi\|_p.
\end{equation}}

The operator $\Xi_{l,m}(\kappa)$ is thus defined through two canonical bilinear forms:
\begin{equation}
 \left\langle \!\left\langle \Xi_{l,m}(\kappa)\varphi, \psi\right\rangle \! \right\rangle = \left\langle \kappa,  \left\langle \! \left\langle \partial_{s_1}^*\cdots \partial_{s_l}^*\partial_{t_1} \cdots \partial_{t_m}\varphi, \psi \right\rangle \! \right\rangle  \right\rangle, \quad \varphi, \psi \in (\cS). 
\end{equation}
This suggest us to employ a formal integral expression:
\begin{align*}
&\Xi_{l,m}(\kappa)
&=\int_{T^{l+m}} \kappa(s_1, \cdots, s_l, t_1,\cdots, t_m) \partial_{s_1}^*\cdots \partial_{s_l}^*\partial_{t_1} \cdots \partial_{t_m}d_{s_1}\cdots d_{s_l}d_{t_1} \cdots d_{t_m}.
\end{align*}

We call $\Xi_{l,m}(\kappa)$ an \textit{integral kernel operator} with kernel distribution $\kappa$.

\ex{Let $(e_n)$ be a complete orthonormal system on $L^2(\R)$. Moreover, let $A\in \cL(\cS(\R), \cS'(\R))$ such that 
$$Ae_k = \sum_{l=1}^{\infty}a_{k,l}e_l.$$
First, we have that
\begin{align*}
    \left\langle \! \left\langle D_{e_k}e_{\mu_{\beta}}(\cdot; \xi), D_{e_k}e_{\mu_{\beta}}(\cdot; \theta) \right\rangle \! \right\rangle &= \left\langle \! \left\langle \langle e_k, \xi \rangle e_{\mu_{\beta}}(\cdot; \xi),  \langle e_k, \theta \rangle e_{\mu_{\beta}}(\cdot; \theta) \right\rangle \! \right\rangle \\
    &= \langle e_k, \xi \rangle \cdot \langle e_k, \theta \rangle \cdot \left\langle \! \left\langle e_{\mu_{\beta}}(\cdot; \xi), e_{\mu_{\beta}}(\cdot; \theta) \right\rangle \! \right\rangle\\
    &= \langle e_k, \xi \rangle \cdot \langle e_k, \theta \rangle \cdot \widehat{I}(\xi, \theta).
\end{align*}
Taking the operator symbol of $\Xi_{1,1}(A) :=\sum_{k=1}^{\infty}\sum_{l=1}^{\infty}D_k^{*}D_{l}a_{k,l}$ yields
\begin{align*}
     \left\langle \! \left\langle \sum_{k=1}^{\infty}\sum_{l=1}^{\infty}D_k^{*}D_{l}a_{k,l}e_{\mu_{\beta}}(\cdot, \xi), e_{\mu_{\beta}}(\cdot, \theta) \right\rangle \! \right\rangle &= \sum_{k=1}^{\infty}\sum_{l=1}^{\infty} a_{k,l}\left\langle \! \left\langle D_k^{*}D_{l}e_{\mu_{\beta}}(\cdot, \xi), e_{\mu_{\beta}}(\cdot, \theta) \right\rangle \! \right\rangle \\
     &=\sum_{k=1}^{\infty}\sum_{l=1}^{\infty} a_{k,l} \left\langle \! \left\langle D_{l}e_{\mu_{\beta}}(\cdot; \xi), D_{k}e_{\mu_{\beta}}(\cdot; \theta) \right\rangle \! \right\rangle \\
     &= \sum_{k=1}^{\infty}\sum_{l=1}^{\infty} a_{k,l}\langle e_l, \xi \rangle \cdot \langle e_k, \theta \rangle \cdot \widehat{I}(\xi, \eta)\\
     &= \sum_{k=1}^{\infty}\sum_{l=1}^{\infty} \langle e_l, Ae_k\rangle \langle e_l, \xi \rangle \cdot \langle e_k, \theta \rangle \cdot \widehat{I}(\xi, \theta) \\
     &= \sum_{k=1}^{\infty}\sum_{l=1}^{\infty} \langle \langle e_l, \xi \rangle e_l,  A\langle e_k, \theta \rangle e_k \rangle \cdot \widehat{I}(\xi, \theta)\\
     &= \langle \sum_{l=1}^{\infty} \langle e_l, \xi \rangle e_l,  \sum_{k=1}^{\infty} A\langle e_k, \theta \rangle e_k \rangle \rangle  \cdot \widehat{I}(\xi, \theta)\\
     &= \langle \xi, A\theta \rangle \cdot \widehat{I}(\xi, \theta). 
\end{align*}}

\section{Acknowledgement}
{\rm{
The financial support provided by the Department of Science and Technology - Accelerated Science and Technology Human Resource Development Program (DOST-ASTHRDP) of the Philippines under the Research Enrichment (Sandwich) Program is gratefully acknowledged. 
}}

\bibliographystyle{alpha}

\end{document}